\documentclass[a4paper,12pt]{amsart}
 
\usepackage{amssymb, amsmath, amsthm, amsfonts}
\usepackage{extsizes}
\usepackage{hyperref}
\usepackage{xcolor,graphicx}

\theoremstyle{plain}
\newtheorem{theorem}{Theorem}[section]
\newtheorem{lemma}[theorem]{Lemma}

\theoremstyle{definition}
\newtheorem{definition}[theorem]{Definition}

\def\r{\mathbb R}
 \def\h{\mathbb H}
 
\usepackage{xcolor,graphicx}

\begin{document}
\title[Homothetical surfaces in hyperbolic space]{Homothetical  
surfaces with constant mean curvature in hyperbolic space}
\author{Rafael Belli}
\address{Department of Mathematics. Federal University of São Carlos. 13565-905 São Carlos, Brazil}
\email{rafaelbelli@estudante.ufscar.br}
\author{Rafael L\'opez}
\address{Department of Geometry and Topology. University of Granada. 18071 Granada, Spain}
\email{rcamino@ugr.es}
\maketitle

\begin{abstract}
We classify all homothetical surfaces with constant mean curvature $H$ in the hyperbolic space $\mathbb{H}^3$. Using the upper half-space model with standard coordinates $(x,y,z)$, these surfaces are defined by the relation $z = \phi(x)\psi(y)$, where $\phi$ and $\psi$ are smooth functions of one variable. We demonstrate that any such surface is necessarily parabolic, meaning that either $\phi$ or $\psi$ is a constant function. Our results cover the minimal case ($H=0$), the case $H^2 \neq 1$, and the critical case $H^2=1$, thereby extending the existing classification of parabolic surfaces in hyperbolic space.
\end{abstract}

\section{Introduction and statement of the main result}

Let $\h^3$ be the $3$-dimensional hyperbolic space, and consider the
upper half-space model of $\h^3$, that is,
$(\mathbb{R}_{+}^{3},g_{\mathbb{H}})$, where $\mathbb{R}
_{+}^{3}=\{(x,y,z)\in \mathbb{R}^{3}:z>0\}$ is endowed with the metric
$$g_{\mathbb{H}} =\frac{1}{z^{2}}g_{\mathbb{E}},$$
 and $g_{\mathbb{E}}=dx^{2}+dy^{2}+dz^{2}$ is the Euclidean metric of $%
\mathbb{R}^{3}$. Using this model, we give the following definition.

\begin{definition} A surface $S$ in $\h^3$ is said to be homothetical if it is locally parametrized by
\begin{equation}\label{ho}
z=\phi(x)\psi(y),
\end{equation}
 where $\phi\colon I\subset\r\to\r$, $\psi\colon J\subset\r\to\r$ are
positive smooth functions of one variable.
 \end{definition}

We emphasize that the definition of a homothetical surface adopted here is model-dependent, as it relies on the Cartesian coordinates of the upper half-space model. However, this choice is motivated by the computational transparency of the metric in these coordinates. The model serves as an optimal workspace where the interplay between Euclidean homotheties and hyperbolic isometries is most evident, facilitating a clearer classification of the surfaces under study. 

This becomes clear with the first examples of
homothetical surfaces, where one of the functions, $\phi$ or $\psi$, is constant. In such cases, $S$ is a parabolic surface in the sense
that it is invariant
 under a one-parameter family of parabolic isometries of $\h^3$. For
example, suppose that $\psi$ is constant and consider the one-parameter group $\mathcal{P}=\{P_t\colon
t\in\r\}$ of parabolic isometries $P_t$, where $P_t(x,y,z)=(x,y+t,z)$. Then $S$ is invariant under $\mathcal{P}$, that is,
 $P_t(S)=S$ for all $t\in\r$. A simple example of a homothetical surface is a
horosphere $z=m$, where $m>0$ is a constant. In this case, we can choose $\phi(x)=1$ and $\psi(y)=m$
for all $x,y$.

\begin{figure}[hbtp]
\begin{center}
\includegraphics[width=.2\textwidth]{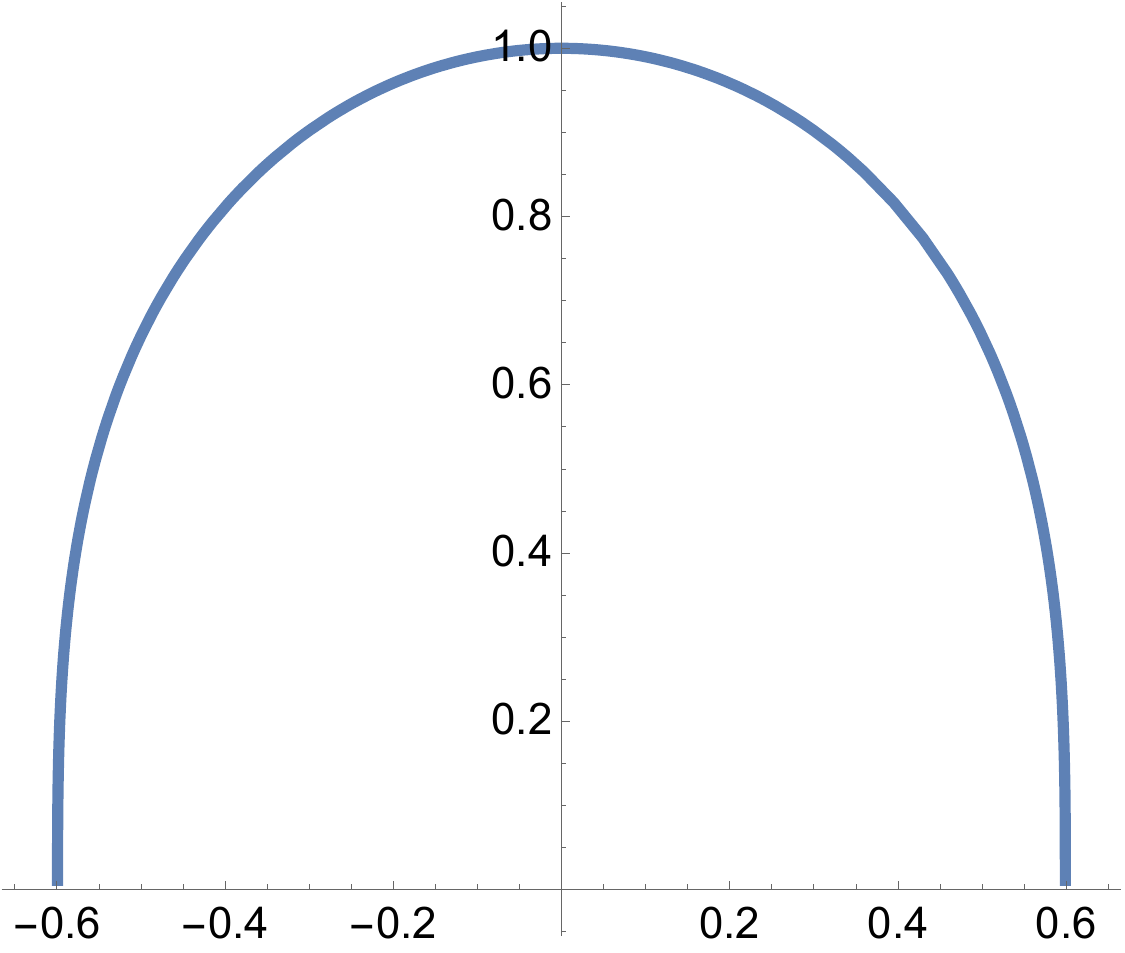}
\includegraphics[width=.25\textwidth]{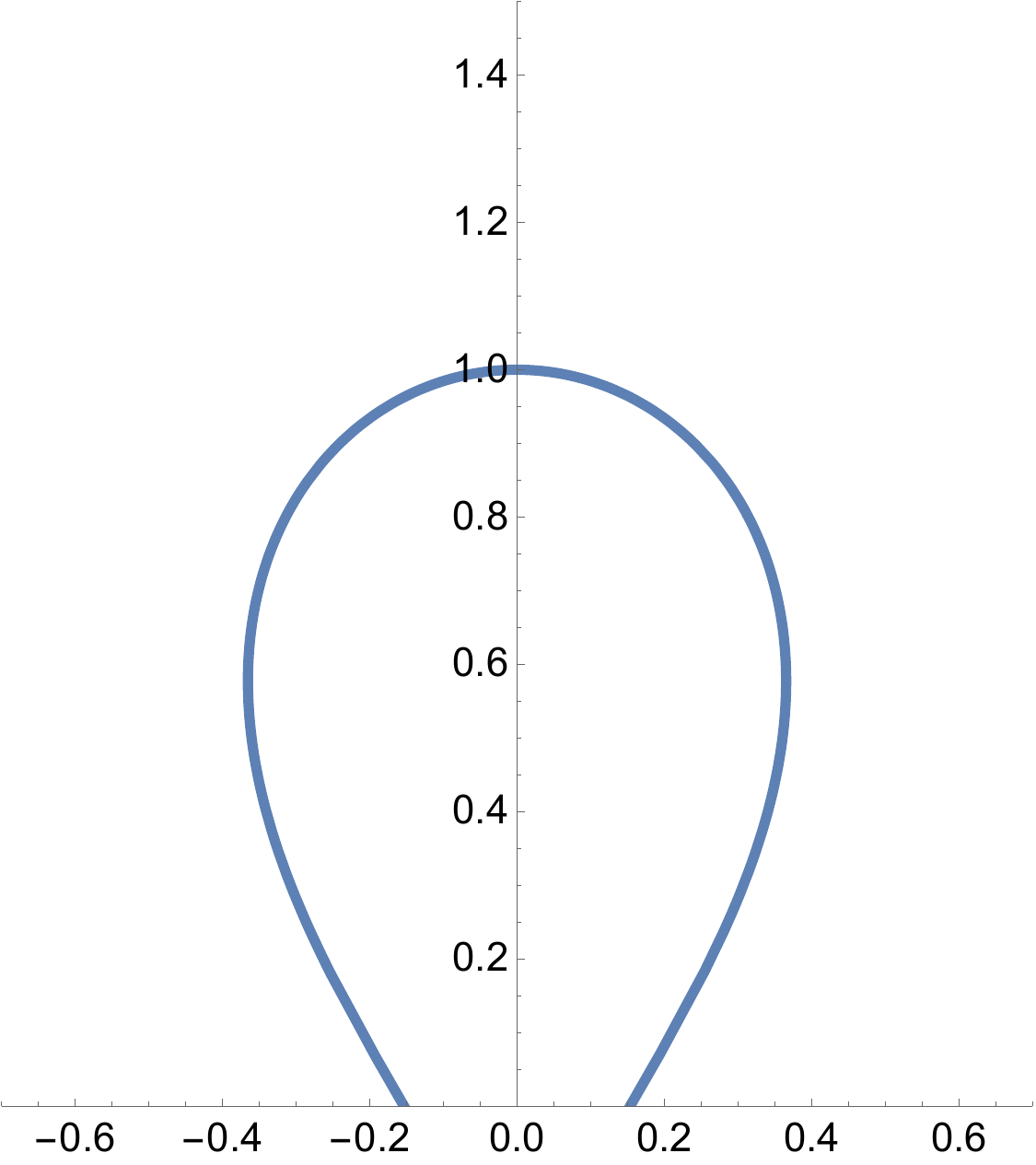}
\includegraphics[width=.27\textwidth]{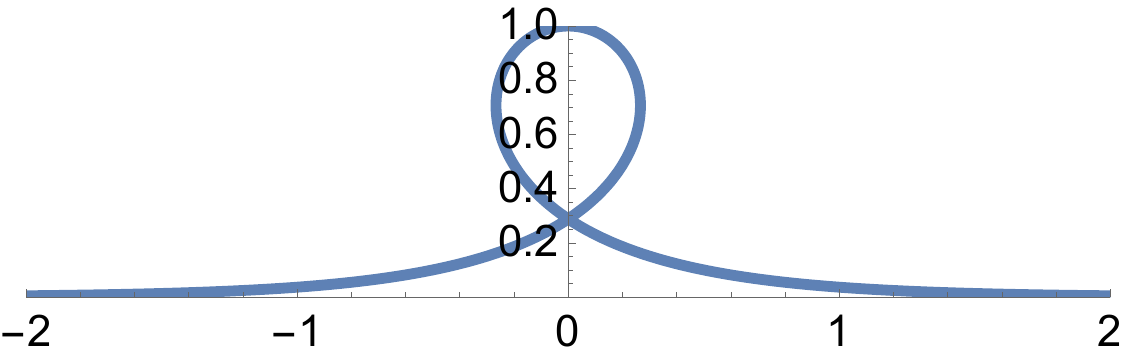}
\includegraphics[width=.22\textwidth]{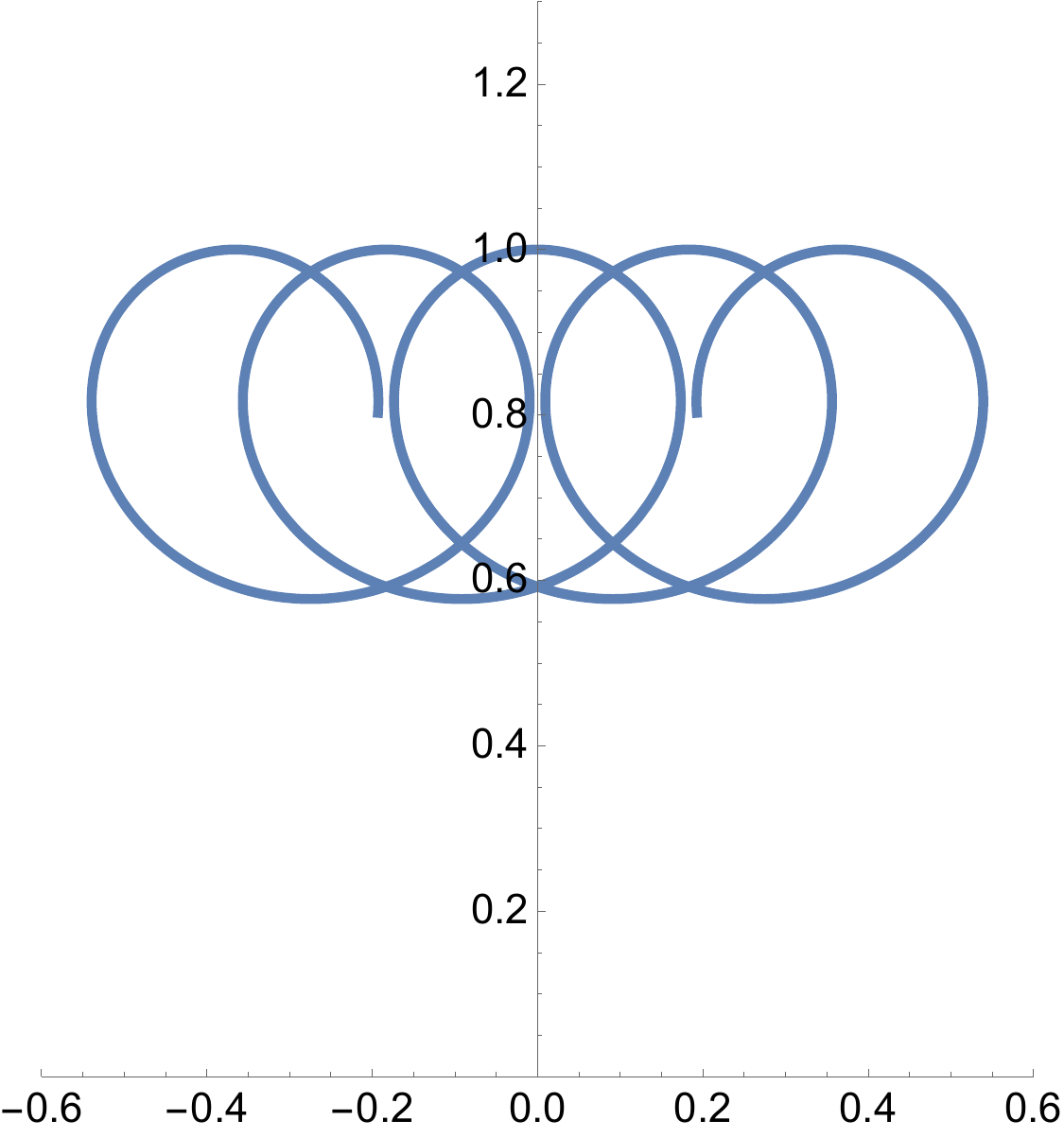}
\end{center}
\caption{Examples of homothetical surfaces of $\h^3$ with constant mean
curvature $H$: from left to right, $H=0$, $H=-\frac12$, $H=-1$ and
$H=2$. } \label{fig1}
\end{figure}

 In this paper, we investigate homothetical surfaces with constant mean
curvature $H$. Note that if $\phi$ or $\psi$ is constant, then the equation $H=c$ reduces to 
a second-order ODE. Indeed, suppose $\psi$ is constant, $\psi(y)=m>0$. Then $S$ is parametrized by
$\Phi(x,y)=(x,y,m\phi(x))$. A straightforward computation shows that the equation $H=c$
can be written as
 $$ \phi \phi'' + 2(1 + \phi'^2)-2H(1 + \phi'^2)^{3/2}=0.$$
 The classification of the parabolic surfaces with constant mean
curvature was carried out by Gomes \cite{go} (see also, \cite{cd,go2}; cf. Fig. \ref{fig1}. Regarding translation surfaces in $\h^3$, that is,
surfaces written as $z=\phi(x)+\psi(y)$, the second author classified
those with zero mean curvature \cite{lo}.

 The purpose of this paper is to provide a full classification of the
homothetical surfaces of $\h^3$ with constant mean curvature.

 \begin{theorem}\label{t1} Let $S$ be a homothetical surface in $\h^3$
parametrized by \eqref{ho}. If $S$ has constant mean curvature, then either 
$\phi$ or $\psi$ is a constant function.
 \end{theorem}

To provide context for this result, we refer to the literature on
homothetical surfaces with constant curvature in other ambient
spaces. The first approach to homothetical surfaces appeared 
in \cite{van}, where minimal, non-degenerate homothetical surfaces in
the Lorentz-Minkowski space were classified. Minimal homothetical
hypersurfaces have been investigated in both Euclidean and semi-Euclidean spaces (\cite{jiu,van 2}). In particular, planes and helicoids are
the only homothetical minimal surfaces of $\r^3$. If we assume that the
mean curvature is constant, the only homothetical surface is the
right-cylinder \cite{ho}. Assuming that the Gauss curvature is constant,
the classification of homothetical surfaces in $\mathbb{R}^{3}$ was completed in \cite{lopez} (see also
\cite{hl0} in the context of separable surfaces in $\r^3$). The study of homothetical surfaces with constant
curvature has been extended to other ambient spaces: \cite{ay1,ay2,mh}.

The proof of Theorem \ref{t1} is divided into different cases according
to the value of $H$. First, in Section \ref{s2}, we obtain a suitable
expression for the mean curvature $H$ for a homothetical surface. Here,
we treat the surface as a special case within the family of separable
surfaces of $\h^3$. The cases studied are: $H=0$ (Section \ref{s3}), 
$H^2\not=0,1$ (Section \ref{s4}), and $H^2=1$ (Section \ref{s5}). In the
computations, we use the software Mathematica. As we
will see, the most difficult case is $H^2=1$. This is to be expected because
the surfaces with constant mean curvature $H=1$ form a rich family of surfaces that 
share similar properties with minimal surfaces in Euclidean space, following the pioneering work of Bryant \cite{br}.

\section{Preliminares}
\label{s2}

In this section, we derive a suitable expression for the mean curvature of a homothetical surface. 
 It is well known that the mean curvature $H$ of a surface in $\mathbb{H}^3$ when   the upper halfspace model is adopted, is related to its Euclidean mean curvature $H_e$ by the formula
\begin{equation} \label{relH}
H = zH_e + N_3,
\end{equation}
where $N=(N_1, N_2, N_3)$ is the Euclidean unit normal vector to the surface \cite[p. 191]{lo0}.

As defined in the introduction, a homothetical surface is given by $z=\phi(x)\psi(y)$. To facilitate the computations, we treat these surfaces as a particular class of separable surfaces. Indeed, by applying the logarithm, we can represent a homothetical surface in the implicit form
\begin{equation*}
f(x)+g(y)+h(z)=0,
\end{equation*}
where
\begin{equation}\label{hh}
\phi(x)=e^{-f(x)},\quad \psi(y) = e^{-g(y)},\quad h(z) = \log(z).
\end{equation}
Note that in this specific case, the function $h$ is fixed as $h(z) = \log(z)$, which implies $z=e^{-f(x)-g(y)}$. Since the claim of Theorem \ref{t1} is that either $\phi$ or $\psi$ is constant, this is equivalent to proving that either $f$ or $g$ is constant. The proof will proceed by contradiction, assuming that neither $f$ nor $g$ is a constant function.

 To generalize this approach, consider three smooth functions $f:I_{1}\to \mathbb{R}$, $g:I_{2}\to \mathbb{R}$, and $h:I_{3}\to \mathbb{R}$. We define a separable surface $S$ as the set
\begin{equation*}
S=\{(x,y,z)\in \mathbb{H}^{3}\colon f(x)+g(y)+h(z)=0\}.
\end{equation*}
If $F(x,y,z)=f(x)+g(y)+h(z)$, the Euclidean unit normal vector field is 
	\[
	N=\frac{\nabla F}{| \nabla F| }=\frac{\left( f^{\prime
		},g',h'\right) }{\sqrt{f'^2+g'^2+h'^2}}
	\]
	
 The Euclidean mean curvature $H_e$ is given by
\begin{equation*}
H_e = -\text{div}\left(\frac{\nabla F}{|\nabla F|}\right) = \frac{\nabla^2 F(\nabla F, \nabla F) - |\nabla F|^2 \Delta F}{2|\nabla F|^3}.
\end{equation*}
Since the Hessian $\nabla^2 F$ is diagonal, the above expression simplifies to
\begin{equation*}
H_e = -\frac{1}{2}\frac{(g'^2+h'^2)f'' + (f'^2+h'^2)g'' + (f'^2+g'^2)h''}{(f'^2+g'^2+h'^2)^{3/2}}.
\end{equation*}
 Substituting these into \eqref{relH}, we obtain the hyperbolic mean curvature:
\begin{equation*}
H = -\frac{z}{2}\frac{(g'^2+h'^2)f'' + (f'^2+h'^2)g'' + (f'^2+g'^2)h''}{(f'^2+g'^2+h'^2)^{3/2}} + \frac{h'}{(f'^2+g'^2+h'^2)^{1/2}}.
\end{equation*}

Let us denote 
\begin{equation*}
u=f(x),\quad v=g(y),\quad w=h(z).
\end{equation*}

From the expression for $h(z)$ in \eqref{hh}, we have 
\begin{equation*}
z=e^{w}.
\end{equation*}
Since $f$ and $g$ are assume to be non-constant, we can introduce new functions of the variables $u$, $v$ and $w$ defined by 
\begin{equation*}
X(u)=f'^{2},\quad Y(v)=g'^{2},\quad Z(w)=h'^{2}=e^{-2w}.
\end{equation*}
In particular, we have 
\begin{equation*}
X'=2f'',\quad Y'=2g'',\quad Z'=2h''=-2 e^{-2w}.
\end{equation*}
The equation for the mean curvature then becomes 
 
\begin{equation} \label{eq1}
(Y+e^{-2w})X'+(X+e^{-2w})Y'-2(X+Y)e^{-2w} -4
e^{-2w}(X+Y+e^{-2w})=-4He^{-w}(X+Y+e^{-2w})^{3/2}.
\end{equation}

Throughout this paper, we will differentiate equations involving functions of $u, v$, and $w$. Since these variables are constrained by $u+v+w=0$, we employ the following lemma:

\begin{lemma}
\label{le1}Let $Q(u,v,w)$ be a smooth function. If $Q(u,v,w)=0$ on the plane $\Pi = \{(u,v,w) : u+v+w=0\}$, then the following holds on $\Pi$:
\begin{equation*}
Q_u = Q_v = Q_w,
\end{equation*}
where $Q_u, Q_v, Q_w$ denote the partial derivatives of $Q$.
\end{lemma}


\section{Case $H=0$}
\label{s3}


Equation \eqref{eq1} can be written as a polynomial equation in $P:=e^{-w}$ given by 
\begin{equation} \label{h0}
A_0(u,v) +A_2(u,v)P^2-4P^4=0,
\end{equation}
where 
\begin{equation}
\begin{split}
A_0&=X'Y+XY', \\
A_2&=-6 (X+Y)+X'+Y'.
\end{split}%
\end{equation}

Note that the leading coefficient associated with $P^4$ is constant. Although the 
coefficients $A_0$ and $A_2$ depend on the variables $u$ and $v$, they actually depend only on on $w=-(u+v)$ owing to Vieta's formulas, which relate the roots of a polynomial to its coefficients.
Viewing Eq. \eqref{h0} as a polynomial in $e^{-2w}$, we then have two roots $%
r_1(w)$ and $r_2(w)$ satisfying 
\begin{equation*}
A_0=-4r_1(w)r_2(w),
\end{equation*}
\begin{equation*}
A_2=4(r_1(w)+r_2(w)).
\end{equation*}
Applying Lemma \ref{le1} to $A_0$ and $%
A_2$, we obtain $(A_0)_u-(A_0)_v=0$ and $(A_2)_u-(A_2)_v=0$, which can be  written as 
\begin{eqnarray}
X''Y-XY''&=&0, \label{h31} \\
X''-6X'&=&Y''-6Y'.
\label{h32}
\end{eqnarray}
From \eqref{h31}, there exists a constant $\lambda\in\mathbb{R}$, such that 
\begin{equation} \label{h33}
\frac{X''}{X}=\lambda=\frac{Y''}{Y}.
\end{equation}

\begin{enumerate}
\item Case $\lambda=0$. From \eqref{h33}, the functions $X$ and $Y$ are linear, but \eqref{h32} implies that the slopes coincide for both. Thus, there exist constants $%
m,a_1,b_1\in\mathbb{R}$ such that 
\begin{equation*}
X(u)=mu+a_1,\quad Y(v)=mv+b_1.
\end{equation*}
Substituting these into \eqref{h0}, we obtain 
\begin{equation*}
m(a_1+b_1-mw)-2\left(3(a_1+b_1)-m-3mw\right)e^{-2w}-4e^{-4w}=0.
\end{equation*}
Since the functions $\{1,w,w e^{-2w},e^{-2w},e^{-4w}\}$ are linearly
independent, we arrive at a contradiction.

\item Case $\lambda=m^2>0$. In this   case,  
\begin{equation*}
X(u)=a_1\cosh(m u)+a_2\sinh(mu),
\end{equation*}
\begin{equation*}
Y(v)=b_1\cosh(mv)+b_2\sinh(mv),
\end{equation*}
with $a_i,b_i\in\r$. Substituting these into \eqref{h32}, we obtain 
\begin{equation}\label{mm}
\begin{split}
0&=(a_2 m-6 a_1) \sinh (m u)+(a_1 m-6 a_2) \cosh (m u)\\
&+(6 b_1-b_2 m) \sinh (m v)+(6 b_2-b_1 m) \cosh (m v).
\end{split}
\end{equation}
It follows that the coefficients must vanish. Since not all $a_i$ and $b_i$ can vanish,
we deduce $m^2=36$, implying $m=\pm 6$. Assuming $m=6$ (a similar argument holds if $m=-6$), from \eqref{mm}, we have $a_2=a_1$ and $b_2=b_1$. Now \eqref{h0}
becomes 
\begin{equation*}
12a_1b_1e^{-6w}-4e^{-4w}=0,
\end{equation*}
which is a contradiction.

\item Case $\lambda=-m^2<0$. In this case,  
\begin{equation*}
X(u)=a_1\cos(m u)+a_2\sin(mu),
\end{equation*}
\begin{equation*}
Y(v)=b_1\cos(mv)+b_2\sin(mv),
\end{equation*}
with $a_i,b_i\in\r$. Substituting into \eqref{h32}, we obtain 
\begin{equation*}
\begin{split}
0&=(-a_2 m+6 a_1) \sin (m u)-(a_1 m+6 a_2) \cos (m u)\\
&+(-6 b_1+b_2 m) \sin (mv)+(6 b_2+b_1 m) \cos (m v).
\end{split}
\end{equation*}
Then, all coefficients must vanish. As in the case $\lambda=m^2$, we deduce that all $a_i$ and $b_i$ are zero, yielding a contradiction.
\end{enumerate}


\section{Case $H\not=0$ and $H^2\not=1$}
\label{s4}


Suppose $H\not=0$. Squaring \eqref{eq1} and collecting terms in $e^{-w}$,
we obtain a polynomial equation of the form

\begin{equation} \label{eq2}
A_0(u,v) +A_2(u,v)e^{-2w}+A_4(u,v)e^{-4w}+A_6(u,v)e^{-6w}+16(1-H^2)e^{-8w}=0,
\end{equation}
where 
\begin{equation}
\begin{split}
A_0&=(Y X'+X Y')^2, \\
A_2&=-2 (24 H^2 (X^2 Y+ X Y^2)+8 H^2 (X^3 + Y^3)+(XY'+X'Y)(6(X+Y)-X'-Y')), \\
A_4&=12((3-4H^2)(X^2+Y^2)+(6-8H^2)XY)\\
&-4(((3X'+5Y')X-(5X'+3Y')Y)+(X'+Y')^2, \\
A_6&=-8 (6(H^2-1) (X+ Y)+X'+Y').
\end{split}%
\end{equation}
Since $H^2\not=1$, the leading coefficient of \eqref{eq2} is constant. Thus, arguing as in the case $H=0$ by using Vieta's relations, the coefficient $A_0$ depends only on $w=-(u+v)$.
Applying Lemma \eqref{le1} to the coefficient $A_0$, we obtain 
\begin{equation*}
X''Y-XY''=0.
\end{equation*}
Thus, there exists $\lambda\in\mathbb{R}$ such that 
\begin{equation*}
\frac{X''}{X}=\lambda=\frac{Y''}{Y}.
\end{equation*}

\begin{enumerate}
\item Case $\lambda=0$. In this case, there exist constants $a_i,b_i$ such that 
\begin{equation*}
X(u)=a_1u+a_2,\quad Y(v)=b_1v+b_2.
\end{equation*}
Applying the lemma again to the coefficient $A_6$, we have 
\begin{equation*}
0=6(H^2-1)(X'-Y')+X''-Y''.
\end{equation*}
Substituting $X$ and $Y$ into this equation, we obtain $a_1=b_1$. Substituting these back into the initial equation %
\eqref{eq2}, we obtain a linear combination of the functions $
w^{3}e^{-2w},w^{2}e^{-2w},w^{2}e^{-4w},w^{2},we^{-2w},we^{-4w},we^{-6w},w,e^{-2w},e^{-4w},\\
e^{-6w},e^{-8w},1
$. Since these functions are linearly independent, all coefficients must vanish. However, the coefficient of $e^{-8w}$ is 
$16(1-H^2)$, which leads to a contradiction.

\item Case $\lambda=m^2>0$. In this case, there exist $a_i,b_i\in\mathbb{R}$ such that 
\begin{equation*}
X(u)=a_1\cosh(mu)+a_2\sinh(mu),
\end{equation*}
\begin{equation*}
Y(v)=b_1\cosh(mv)+b_2\sinh(mv).
\end{equation*}
We now consider the coefficient $A_2$ of \eqref{eq2} and  substitute $X$
and $Y$. Applying the lemma, we obtain a polynomial equation on 
$$%
\{\cosh(mu),\sinh(mu),\cosh(mv),\sinh(mv)\},$$
 with $0\leq m\leq 3$. In order
to simplify the computations and the resulting expressions, we set $v=0$. This yields 
$ \sum_{n=0}^3A_n\cosh(nu)+B_n\sinh(nu) =0$. 
All coefficients must vanish. From $A_3$ and $B_3$, we obtain 
\begin{equation*}
\begin{split}
A_3&=6a_2mH^2(3a_1^2+a_2^2),\\
 B_3&=6a_1mH^2(a_1^2+3a_2^2).
 \end{split}
\end{equation*}
Thus $a_1=a_2=0$, which implies $X=0$, a contradiction.

\item Case $\lambda=m^2<0$. Similarly, there exist $a_i,b_i\in\mathbb{R}$ such that 
\begin{equation*}
X(u)=a_1\cos(mu)+a_2\sin(mu),
\end{equation*}
\begin{equation*}
Y(v)=b_1\cos(mv)+b_2\sin(mv).
\end{equation*}
The argument is similar to the case $\lambda=-m^2$. By considering the coefficient $%
A_2$ of \eqref{eq2}, and setting $v=0$, we arrive at the equation 
$ \sum_{n=0}^3A_n\cos(nu)+B_n\sin(nu) =0$. Then, 
\begin{equation*}
\begin{split}
A_3&=6a_2mH^2(3a_1^2-a_2^2),\\
 B_3&=6a_1mH^2(-a_1^2+3a_2^2).
 \end{split}
\end{equation*}
Since $A_3$ and $B_3$ must vanish, we conclude that $a_1=a_2=0$ again, which 
is a contradiction.
\end{enumerate}


\section{Case $H^2=1$}
\label{s5}


Equation \eqref{eq2}   reduces to
\begin{equation}
A_{0}(u,v)+A_{2}(u,v)e^{-2w}+A_{4}(u,v)e^{-4w}+A_{6}(u,v)e^{-6w}=0,
\label{eq21}
\end{equation}
where

\begin{equation}
\begin{split}
A_{0}& =-( YX'+XY') ^{2}, \\
A_{2}& =48(XY^2+X^2Y)+16(X^3+Y^3)+2(XY'+YX')(6(X+Y)-X'-Y'), \\
A_{4}& =12(XX'+YY')+20(XY'+X'Y)-(
X'+Y') ^{2}+12(X+Y)^{2}, \\
A_{6}& =8(X'+Y').
\end{split}%
\end{equation}

\textit{Claim:} $X'+Y'\not=0$.

The proof is by contradiction. Suppose $X''=0$.
Then there exists $a\in\mathbb{R}$ such that $X'=a=-Y'$. Thus $%
X(u)=au+b$ and $Y(u)=-av+c$, for some $b,c\in\mathbb{R}$. Substituting these into \eqref{eq21},
and next, multiplying by $e^{4w}$, we arrive at 
\begin{equation*}
B_0+B_2e^{2w}+B_4 e^{4w}=0,
\end{equation*}
where 
\begin{equation*}
\begin{split}
B_0&= 4 (a^2 (3 u^2-2 u (3 v+1)+v (3 v-2))+2 a (-3 v (b+c)+3
b u-b+3 c u+c)+3 (b+c)^2), \\
B_2&=4 (a (u-v)+b+c) \Big((a^2 (4 u^2-u (8 v+3)+v (4 v-3)) \\
&+a (-8 v (b+c)+8 b u-3 b+8 c u+3 c)+4 (b+c)^2\Big), \\
B_4&= -a^2 (a (u+v)+b-c)^2= -a^2 (a w+b-c)^2.
\end{split}%
\end{equation*}
Note that $B_4\not=0$. Consequently, $\frac{B_0}{B_4}$ and $\frac{B_2}{B_4}$ depend only
 on $u+v$. Applying Lemma \ref{le1} to $\frac{B_0}{B_4}$, we have 
\begin{equation*}
(B_0)_uB_4-B_0(B_4)_u-(B_0)_vB_4+B_0(B_4)_v=0.
\end{equation*}
This is a polynomial equation in $u$ of the form $\sum_{n=0}^3C_n(v)u^n=0$,
where $C_3=48a^6$. Thus $a=0$, which implies $X(u)=b$ and $Y(v)=c$. In this case 
\eqref{eq21} becomes $12(b+c)^2+16(b+c)^3e^{-2w}=0$, which implies $c=-b$. It follows that 
\begin{equation*}
f'^{2}=b=-g'^{2},
\end{equation*}
which leads to $b=0$. Thus $X=Y=0$, which is a contradiction.

Having established the claim, we can divide by $A_6$ in the equation %
\eqref{eq21}, obtaining 
\begin{equation} \label{eq211}
Q_0(u,v) +Q_2(u,v)e^{-2w}+Q_4(u,v)e^{-4w}+ e^{-6w} =0,
\end{equation}
where 
\begin{equation*}
Q_0=\frac{A_0}{A_6},\quad Q_2=\frac{A_2}{A_6},\quad Q_4=\frac{A_4}{A_6}.
\end{equation*}
Since the leading coefficient is constant, it follows from Vieta's formulas that the
coefficients $Q_0$, $Q_2$ and $Q_4$ must depend only on the roots $r_{1}(
w) $, $r_{2}(w) $ and $r_{3}(w) $ of %
\eqref{eq211}. In particular, we have 
\begin{equation*}
Q _{0}(u+v) = \frac{-( XY'+YX') ^{2}}{8(X'+Y')} .
\end{equation*}
By the lemma, we have $(Q_0)_u-(Q_0)_v=0$ which, after
simplification, can be written as 
\begin{equation} \label{eq200}
(Y X'+X Y') (Y''(X'(Y-2
X)-X Y')+X''(Y X'-Y'(X-2
Y)))=0.
\end{equation}

\subsection{Case $Y X'+X Y'=0$ identically}

In this case, there exists $a\in 
\mathbb{R}$ such that 
\begin{equation*}
\frac{X'}{X}=a=-\frac{Y'}{Y}.
\end{equation*}
The case $a=0$ is precluded because then $X'+Y'=0$,
which was previously discarded. Thus $a\not=0$, yielding 
\begin{equation*}
X(u)=e^{a u},\quad Y(v)=e^{-a v}.
\end{equation*}
Furthermore, we also have $(Q_2)_u-(Q_2)_v=0$. This equation becomes 
\begin{equation*}
e^{au}-e^{-av}=0,
\end{equation*}
yielding a contradiction.

 \subsection{Case $YX'+XY'\not=0$ on some open set $U_{0}\times
V_{0}$ of $I\times J$.}
From \eqref{eq200}, we have 
\begin{equation}
Y''( X'(Y-2X)-XY') +X'' ( YX'-Y'(X-2Y)) =0. \label{eq3}
\end{equation}
Dividing by $XX'YY'$, we obtain 
\begin{equation*}
\frac{2X''}{XX'}-2\frac{Y''}{
YY'}-\frac{X''}{X'Y}+\frac{X''}{XY'}-\frac{Y''}{X'Y}+\frac{Y''} {XY'}=0.
\end{equation*}

Differentiating with respect to $u$ and then with respect to $v$ yields 
\begin{equation*}
-\left(\frac{X''}{X'}\right)'\left(\frac{1}{Y}\right)'+\left( 
\frac{X''}{X}\right)'\left(\frac{1}{Y'}\right)'-\left( 
\frac{1}{X'}\right)'\left(\frac{Y''}{Y}\right)'+\left( 
\frac{1}{X}\right)'\left(\frac{Y''}{Y'}\right)'=0,
\end{equation*}
which can be rewritten as 
\begin{equation*}
-\frac{\left(\frac{X''}{X'}\right)'}{\left(\frac{1}{
X'}\right)'}\frac{\left(\frac{1}{Y}\right)'}{\left(\frac{1}{Y'}
\right)'}+\frac{\left(\frac{X''}{X}\right)'}{\left(\frac{1}{
X'}\right)'}-\frac{\left(\frac{Y''}{Y}\right)'}{\left( 
\frac{1}{Y'}\right)'}+\frac{\left(\frac{Y''}{Y'}\right)'}{\left(\frac{1}{Y'}\right)'}\frac{\left(\frac{1}{X}\right)'}{\left(\frac{1}{X'}\right)'}=0.
\end{equation*}

Differentiating once more with respect to $u$ and $v$, we find
\begin{equation*}
\dfrac{\left( \frac{\left(\dfrac{X''}{X'}\right)'}{\left(\dfrac{ 1}{X'}\right)'}\right)'}{\left(\frac{\left(\dfrac{1}{X}\right)'}{\left( \dfrac{ 1}{X'}\right)'}\right)'}
=\dfrac{\left( \frac{\left(\dfrac{Y''}{Y'}\right)'}{\left(\dfrac{ 1}{Y'}\right)'}\right)'}{\left(\frac{\left(\dfrac{1}{Y}\right)'}{\left( \dfrac{ 1}{Y'}\right)'}\right)'}.
\end{equation*}
Therefore, there exist constants $\lambda,a_1,b_1\in\mathbb{R}$ such that 
\begin{equation*}
\frac{\left(\dfrac{X''}{X'}\right)'}{\left(\dfrac{1}{
X'}\right)'} =\lambda\frac{\left(\dfrac{1}{X}\right)'}{\left(\dfrac{1}{
X'}\right)'}+a_1,\qquad \frac{\left(\dfrac{Y''}{
Y'}\right)'}{\left(\dfrac{1}{Y'}\right)'} =\lambda\frac{\left( 
\dfrac{1}{Y}\right)'}{\left(\dfrac{1}{Y'}\right)'}+b_1,
\end{equation*}
or, equivalently, 
\begin{equation*}
\begin{split}
\left(\dfrac{X''}{X'}\right)'&=\lambda\left(\dfrac{1}{X}
\right)'+a_1\left(\dfrac{1}{X'}\right)',\\
 \left(\dfrac{Y''} {Y'}\right)'&=\lambda\left(\dfrac{1}{Y}\right)'+b_1\left( 
\dfrac{1 }{Y'}\right)'.
\end{split}
\end{equation*}
Integrating, we find that there are $a_2,b_2\in\mathbb{R}$ such that 
\begin{equation} \label{xy}
\begin{split}
X''&=\lambda \frac{X'}{X}+a_2X'+a_1 \\
Y''&=\lambda \frac{Y'}{Y} +b_2 Y'+b_1.
\end{split}%
\end{equation}

\begin{enumerate}
\item Case $\lambda =0$. In this case, there exist $a_3 ,a_4,b_3 ,b_4 \in \mathbb{R} $,
with $a_3 b_3 \not=0$, such that 
\begin{equation*}
X(u)=a_3 e^{a_2 u}-\frac{a_1}{a_2} u+a_4,\qquad Y(v)=b_3 e^{b_2 v}- \frac{
b_1 }{b_2} v+b_4 .
\end{equation*}
Note that the roles of $u$ and $v$ are symmetric. Substituting into \eqref{eq3},
we obtain the following: 
\begin{equation}
S_{1} ue^{ a_2 u} +S_{2} e^{ a_2 u} +S_{3} e^{2 a_2 u} +S_{4} u+S_{5}=0,
\label{nova 1}
\end{equation}
where the coefficients are linear combinations of the functions $\{1,v,
e^{b_2 v}, v e^{b2_2 v},e^{2b_2 v}\}$. Since $\left\{ e^{ a_2 u} u,e^{ a_2
u},e^{2 a_2 u},u,1\right\} $ is a linearly independent set, all coefficients must be all zero. Computing $S_3$ gives 
\begin{equation*}
S_{3} = ( -a_2 ^{3}a_3 ^{2}b_2 ( b_3 a_2 b_2 ^{2}-b_3 a_2 ^{2}b_2
+2b_3 b_2 ^{3}) ) e^{ b_2 v}+( -a_2 ^{5}a_3 ^{2}b_1 b_2
) v+a_2 ^{3}a_3 ^{2}b_2 ( b_2 b_4 a_2 ^{2}+b_1 a_2 ).
\end{equation*}
Since $S_3=0$, we obtain 
\begin{equation*}
( -a_2 ^{2}b_1 ) v+( b_3 a_2 ^{2}b_2 -b_3 a_2 b_2 ^{2}-2b_3 b_2 ^{3}) e^{
b_2 v} +( b_2 b_4 a_2 ^{2}+b_1 a_2 ) =0.
\end{equation*}
Using the linear independence of $\{ 1, v,e^{ b_2 v} \} $, it follows that 
\begin{equation*}
\begin{split}
0&=b_1 , \\
0&= ( a_2 +b_2 ) ( a_2 -2b_2 ), \\
0&= b_1 +a_2 b_2 b_4,
\end{split}%
\end{equation*}
which implies $ b_1=b_4 =0 $ and $a_2 =\epsilon b_2 $, where either 
$\epsilon =-1$ or $\epsilon =2$. By symmetry, we deduce $a_1=a_4=0$ and $b_2=\epsilon a_2$. If $\epsilon=2$, then $a_2=b_2$, which is a contradiction. Thus $a_2=-b_2$, which leads to 
\begin{equation*}
X(u)=a_3 e^{-a_2 u},\qquad Y(v)=b_3 e^{b_2 v}.
\end{equation*}
In such a case, it is immediate that $X'Y+XY'=0$,
another contradiction.

\item Case $\lambda \not=0$. A first situation that we will discard is that $X$ and $Y$ cannot be linear. This will be crucial in the rest of arguments. We emphasize this case in the following claim:

{\it Claim L:} The functions $X$ and $Y$ cannot be linear functions.

The proof is the following. Suppose $X(u)=au+b$, $a,b\in\r$. Then \eqref{xy} implies
\[ 0=a( a_1 +aa_2 ) u+( a\lambda +b(a_1 +aa_2) ).\]
This implies $a_1 +aa_2=0$ and $\lambda=0$, contradicting our assumption $\lambda\not=0$. A similar argument proves that $Y$ cannot be a linear function.

We now begin with the proof of the case $\lambda\not=0$. Applying Lemma \ref{le1} in $Q_0$, we have $(Q_0)_u-(Q_0)_v=0$. Substituting $%
X''$ from \eqref{xy}, we find 
\begin{equation}
\begin{split}
 0 &=\lambda Y X'^{2}-2X^{2}X'Y''-X^{2}Y'Y''-a_{1}X^{2}Y'-a_{2}X^{2}X'Y'+XYX'Y''+a_{1}XYX'
 \\
&+2a_{1}XYY'-\lambda X X'Y'+2\lambda Y X'Y'+a_{2}XYX'^{2}+2a_{2}XYX'Y',
\label{nova 4}
\end{split}
\end{equation}

\end{enumerate}
 
For a fixed $v$, (\ref{nova 4}) can be viewed as a family of two variable polynomial equations
with respect to $X$ and $X'$, that is,%
\begin{equation}
\begin{split}
 0 &=p(X,X')\\
&=(p_{3}X+p_{6})X'^{2}+(p_{1}X^{2}+p_{4}X+p_{7})X'+(p_{2}X^{2}+p_{5}X) \label{nova 5}
\end{split}
\end{equation}%
where the coefficients $p_i(v)$ are given by 
\begin{equation}\label{ps}
\begin{split}
p_{1}(v) &=-2Y''-a_{2}Y', \\
p_{2}(v) &=-Y'Y''-a_{1}Y', \\
p_{3}(v) &=a_{2}Y, \\
p_{4}(v) &=YY''-\lambda Y'+2a_{2}YY'+a_{1}Y, \\
p_{5}(v) &=2a_{1}YY', \\
p_{6}(v) &=\lambda Y, \\
p_{7}(v) &=2\lambda YY'.
\end{split}
\end{equation}

We   utilize a result on the intersections of zeros of two polynomial. For $%
P \in \mathbb{R}[X,Y]$, let $\mathcal{Z}(P)$ denote its set of zeros. The following result can be found in \cite{fultonnn}:

\begin{lemma}
\label{fulton} If $P,Q\in \mathbb{R}[X,Y]$ have no common factors, then $%
\mathcal{Z}(P)\cap \mathcal{Z}(Q)$ is a finite set.
\end{lemma}

We analyze the possible decompositions of $p(X,X')$ in \eqref{nova 5}. For each possible decomposition, we obtain at least one differential equation depending only on $Y$. These equations are the so-called factorability conditions. We will obtain four equations $E_n$, $1\leq n\leq 4$. With these equations at hand, we proceed with the following strategy: first, assume that $E_1$ holds throughout an open interval $V_1$ in the $v$ variable. After resolving this case, we assume that $E_1$ does not hold on all of $V_1$; then there exists a neighborhood $V_2\subset V_1$ on which $E_1$ fails. We then assume that $E_2$ holds throughout $V_2$, and so on. At the end of this procedure, we find a neighborhood $V_4$ on which $p$ must be irreducible for all $v \in V_4$. We will use Lemma \ref{fulton} to solve this case.\\\\
First, if $p_1=0$ in some neighborhood of $V_0$, then by \eqref{xy} we have $Y'=-\frac{a_2}{2}Y$. Using this fact in \eqref{eq3}, we conclude that
$X'=\frac{a_2}{2}X$. However, this implies that $XY' + X'Y = 0$, which is a contradiction. Moreover, if $p_2=0$ in some neighborhood of $V_0$, then by \eqref{xy} we obtain $Y' = 0$, which is also a contradiction.

Therefore, we may assume that there exists some neighborhood $V_1 \subset V_0$ on which $p_1 p_2 p_3 \neq 0$ because $p_3\not=0$. Recall that $V_0$ is an interval where $YX'+XY'\not=0$ in $U_0\times V_0\subset I\times J$. Thus, there will be no change in the degree of each polynomial with respect to $X$ in \eqref{nova 5}.

Since $\deg(p(X,X'))=2$ in \eqref{nova 5} with respect to $X'$, there are two possible decompositions of $p$:
\begin{eqnarray}
p(X,X') &=&j( q_{2}X'^{2}+q_{1}X'+q_{0}), \label{de1}\\
p(X,X') &=&(m_{1}X'+m_{0})(n_{1}X'+n_{0}),\label{de2}
\end{eqnarray}
where $q_0$, $q_1$, $q_2$, $m_0$, $m_1$, $n_0$, $n_1$ and $j$ are polynomials with respect to $X$.
\begin{enumerate}
\item Assume the first decomposition \eqref{de1}. By (\ref{nova 5}), 
we obtain the following equations
\begin{eqnarray}
 jq_{2} &=&p_{3}X+p_{6}, \label{dec 1 eq 3} \\
jq_{1} &=&p_{1}X^{2}+p_{4}X+p_{7}, \label{dec 2 eq 3}\\
jq_{0} &=&p_{2}X^{2}+p_{5}X. \label{dec 3 eq 3}
\end{eqnarray}
By (\ref{dec 1 eq 3}), we have
\begin{equation*}
j=\frac{p_{3}}{q_{2}}X+\frac{p_{6}}{q_{2}}.
\end{equation*}%
Moreover, we conclude that $\deg (q_{1})=\deg
(q_{0})=1$, that is:
\begin{equation*}
\begin{split}
q_{1}&=q_{1,1}X+q_{1,0},\\
 q_{0}&=q_{0,1}X+q_{0,0}.
 \end{split}
\end{equation*}%
Using the expression of $q_1$ together with $j$ in (\ref{dec 2 eq 3}), we find
\begin{eqnarray*}
0 &=&jq_{1}-(p_{1}X^{2}+p_{4}X+p_{7}) \\
&=&\frac{1}{q_{2}}(p_{3}X+p_{6})(q_{1,1}X+q_{1,0})-(p_{1}X^{2}+p_{4}X+p_{7})
\\
&=&\frac{1}{q_{2}}(p_{3}q_{1,1}-p_{1}q_{2})X^{2}+\frac{1}{q_{2}}%
(p_{3}q_{1,0}-p_{4}q_{2}+p_{6}q_{1,1})X-\frac{1}{q_{2}}%
(p_{7}q_{2}-p_{6}q_{1,0}).
\end{eqnarray*}

The coefficients must all be zero, which implies 
\begin{equation}\label{y1}
0=p_{1}p_{6}^{2}+p_{3}^{2}p_{7}-p_{3}p_{4}p_{6}=(a_{2}Y+2\lambda )Y''+a_{1}a_{2}Y.
\end{equation}

\item Assume the second decomposition \eqref{de2}. Then
\[ p(X,X') =m_{1}n_{1}X'^{2}+(m_{0}n_{1}+m_{1}n_{0})X'+m_{0}n_{0}.\] 
By (\ref{nova
5}), we have the following equations 
\begin{eqnarray}
m_{1}n_{1} &=&p_{3}X+p_{6}, \label{dec 1 eq 4} \\
m_{0}n_{1}+m_{1}n_{0} &=&p_{1}X^{2}+p_{4}X+p_{7}, \label{dec 2 eq 4} \\
m_{0}n_{0} &=&p_{2}X^{2}+p_{5}X. 
\label{dec 3 eq 4}
\end{eqnarray}

By (\ref{dec 1 eq 4}), $\deg(m_1)+\deg(n_1)=1$. By symmetry we can assume without loss of generality that $\deg (m_{1})=1$ and $%
\deg (n_{1})=0$. Finally, by (\ref{dec 3 eq 4}), $\deg (m_{0})+\deg (n_{0})=2$. Using these facts in (\ref{dec 2 eq 4}) we conclude that the only possibilities for $m_0$ and $n_0$ are either $(\deg (m_{0}),\deg (n_{0}))=(2,0)$
or $(\deg (m_{0}),\deg (n_{0}))=(1,1)$.
\begin{enumerate}
\item If $(\deg (m_{0}),\deg (n_{0}))=(2,0)$, the expressions are 
\begin{equation*}
\begin{split}
m_{0}& =m_{0,2}X^{2}+m_{0,1}X+m_{0,0}, \\
m_{1}& =\frac{p_{3}}{n_{1}}X+\frac{p_{6}}{n_{1}}.
\end{split}%
\end{equation*}

Using these expressions in (\ref{dec 2 eq 4}), we find
\begin{eqnarray*}
 0&=&m_{0}n_{1}+m_{1}n_{0}-( p_{1}X^{2}+p_{4}X+p_{7})\\ 
&=&( \frac{p_{2}}{n_{0}}X^{2}+\frac{p_{5}}{n_{0}}X) n_{1}+\frac{%
 p_{3}X+p_{6}}{n_{1}}n_{0}-( p_{1}X^{2}+p_{4}X+p_{7})\\ 
&=& ( \frac{1}{n_{0}}n_{1}p_{2}-p_{1}) X^{2}+( 
\frac{n_{0}}{n_{1}}p_{3}-p_{4}+\frac{1}{n_{0}}n_{1}p_{5}) 
X+( \frac{n_{0}}{n_{1}}p_{6}-p_{7})
\end{eqnarray*}
That is, we have a polynomial equation with respect to $X$. All coefficients must be zero. By the first we have $\frac{n_1}{n_0}=\frac{p_1}{p_2}$. Using it in the last coefficient we have 
\begin{eqnarray*}
0=p_1p_7-p_2p_6,
\end{eqnarray*}
which implies the following equality:
\begin{equation}\label{y2}
0=3Y''+2a_{2}Y'-a_{1}.
\end{equation}

\item If $(\deg (m_{0}),\deg (n_{0}))=(1,1)$, by (\ref{dec 1 eq 4}) and (\ref{dec 3 eq 4}) we have the following
expressions
\begin{equation*}
\begin{split}
m_{1}& =\frac{p_{3}}{n_{1}}X+\frac{p_{6}}{n_{1}}, \\
m_{0}& =m_{0,1}X+m_{0,0}, \\
n_{0}& =n_{0,1}X+n_{0,0}.
\end{split}%
\end{equation*}

By (\ref{dec 3 eq 4}), we have the following
\begin{eqnarray}
0 &=&m_{0}n_{0}-( p_{2}X^{2}+p_{5}X) \label{vinte oito}\\
&=&( m_{0,1}X+m_{0,0}) ( n_{0,1}X+n_{0,0}) -(
p_{2}X^{2}+p_{5}X)\notag \\
&=&\allowbreak ( m_{0,1}n_{0,1}-p_{2}) X^{2}+(
m_{0,0}n_{0,1}-p_{5}+m_{0,1}n_{0,0}) \allowbreak X+m_{0,0}n_{0,0}\notag
\end{eqnarray}
Again we have a polynomial equation with respect to $X$. Since all coefficients must be zero, by the last we have two possibilities: Either $m_{0,0}=0$ or $n_{0,0}=0$.\\
If $m_{0,0}\not =0$ and $n_{0,0}=0$, \eqref{vinte oito} becomes
\begin{eqnarray*}
 0= (m_{0,1}n_{0,1}-p_2)X^2+(m_{0,0}n_{0,1}-p_5) X,
\end{eqnarray*}
which implies $m_{0,1}=\frac{p_2}{n_{0,1}}$ and $m_{0,0}=\frac{p_5}{n_{0,1}}$. Using these equalities in (\ref{dec 2 eq 4}) we find
\begin{eqnarray*}
0 &=&m_{0}n_{1}+m_{1}n_{0}-( p_{1}X^{2}+p_{4}X+p_{7}) \\
&=&( \frac{p_{2}}{n_{0,1}}X+\frac{p_{5}}{n_{0,1}}) n_{1}+( 
\frac{p_{3}}{n_{1}}X+\frac{p_{6}}{n_{1}}) (
n_{0,1}X+n_{0,0}) -( p_{1}X^{2}+p_{4}X+p_{7}).
\end{eqnarray*}
That is, we the the following polynomial equation with respect to $X$:
\begin{eqnarray*}
 0=-n_{0,1}( n_{1}p_{1}-p_{3}n_{0,1}) X^{2}+(
p_{2}n_{1}^{2}-p_{4}n_{1}n_{0,1}+p_{6}n_{0,1}^{2}) \allowbreak
X+( p_{5}n_{1}^{2}-p_{7}n_{0,1}n_{1}).
\end{eqnarray*}
Since all coefficients must be zero, by the first and last coefficients we find the following equation:
\begin{equation*}
p_{1}p_{7}-p_{3}p_{5}=0,
\end{equation*}
which implies the following equation:
\begin{equation*}
 2\lambda Y''+\lambda a_2Y'+a_1a_2Y=0.
\end{equation*}
On the other hand, if $m_{0,0}=0$, \eqref{vinte oito} turns into
\begin{eqnarray*}
0 &=&( m_{0,1}n_{0,1}-p_{2}) X^{2}+(
-p_{5}+m_{0,1}n_{0,0})X ,
\end{eqnarray*}
since all coefficients must be zero, we have $n_{0,1}=\frac{p_{2}}{m_{0,1}}$ and $n_{0,0}=\frac{p_{5}}{m_{0,1}}$. Using these equalities again in (\ref{dec 2 eq 4}) we find
\begin{eqnarray*}
0 &=&( m_{0}n_{1}+m_{1}n_{0}) -(
p_{1}X^{2}+p_{4}X+p_{7}) \\
&=&( ( m_{0,1}X) n_{1}+( \frac{p_{3}}{n_{1}}X+\frac{%
p_{6}}{n_{1}}) ( \frac{p_{2}}{m_{0,1}}X+\frac{p_{5}}{m_{0,1}}%
) ) -( p_{1}X^{2}+p_{4}X+p_{7}).
\end{eqnarray*}
That is, we have the following polynomial equation with respect to $X$:
\begin{equation*}
 0=( p_{2}p_{3}-n_{1}p_{1}m_{0,1}) X^{2}+(
n_{1}^{2}m_{0,1}^{2}-p_{4}n_{1}m_{0,1}+p_{2}p_{6}+p_{3}p_{5})
\allowbreak X+( p_{5}p_{6}-n_{1}p_{7}m_{0,1}).
\end{equation*}
Since all coefficients are zero, by the first we find $m_{0,1}=\frac{p_{5}p_{6}}{n_{1}p_{7}}=Y\frac{a_{1}%
}{n_{1}}$ (in particular, $a_{1}\neq 0$). Using this fact in the last coefficient we find
\begin{equation*}
 p_{1}p_{5}p_{6}-p_{2}p_{3}p_{7}=0,
\end{equation*}
which implies the following equation:
\begin{equation*}
 -2a_{1}+Y'a_{2}=0.
\end{equation*}
In conclusion, for the second possible decomposition of $p$ one of the following equations must hold:
\begin{equation}\label{y3}
0=2\lambda Y''+\lambda a_{2}Y'+a_{1}a_{2}Y,
\end{equation}
or
\begin{equation}\label{y4}
 0=a_{2}Y'-2a_{1},
\end{equation}
where the first equality holds if $m_{0,0}\neq 0$ and the second if $%
m_{0,0}=0$. Note that $a_{1}\neq 0$ in the second equation, because $\deg
m_{0}=1$ and $m_{0,1}=\frac{a_{1}}{n_{1}}Y$.
\end{enumerate}
\end{enumerate}

Consequently, from \eqref{y1}, \eqref{y2}, \eqref{y3} and \eqref{y4}, we have four possibilities. 
\begin{description}
\item[Case 1] Suppose $(a_{2}Y+2\lambda )Y''+a_{1}a_{2}Y=0$ on all of $V_1$. Then we can isolate $Y''$ almost everywhere in $V_{1}$, yielding
\begin{equation}
Y''=\frac{-a_{1}a_{2}Y}{a_{2}Y+2\lambda }. \label{y'' expression}
\end{equation}%
Combining this with (\ref{xy}), we conclude 
\begin{equation*}
Y'=-\frac{a_{2}(a_{1}+b_{1})Y^{2}+2\lambda b_{1}Y}{(b_{2}Y+\lambda
)(a_{2}Y+2\lambda )}
\end{equation*}
almost everywhere in $V_{1}$. Substituting these values of $Y''$ and $Y'$ in \eqref{eq3}, along with $X''$ from \eqref{xy}, we obtain an expression of the type $\sum_{n=0}^3 C_n(u) Y^n=0$. Thus, the coefficients $C_n$ must vanish. The coefficients $C_n$ are polynomials on $X$ and $X'$. The equation $C_3=0$ is 
$$b_2(\lambda+a_2X)X'^2-2(a_1+b_1)(\lambda+a_2X)X'-2a_1(a_1+b_1)X=0.$$
First of all, suppose $(b_2,a_1+b_1)\neq (0,0)$. Solving for $X$ and substituting into $C_0=0$ and $C_2=0$, we obtain two polynomial equations in $X'$, namely, 
\begin{equation}\label{di1}
\begin{split}
0&=8(a_1+b_1)(2a_1a_2+2a_2b_1+a_1b_2-b_1b_2) \\
&+4b_2(-2a_1a_2+b_1(b_2-2a_2)) X',\\
0&=-4b_1(a_1+b_1)(a_1a_2+a_2b_1+2a_1b_2)\\
&+2(a_1+b_1)(2a_2+b_2)(a_1a_2+a_2b_1+2a_1b_2)X'\\
&-2a_2b_2(a_1a_2+a_2b_1+2a_1b_2)X'^2=0.
\end{split}
\end{equation}
All coefficients of $X'$ must vanish identically. Since $a_1a_2\not=0$, if $b_1=-a_1$, \eqref{di1} leads to $-4a_1b_2^2 X'=0$, implying $b_2=0$, a contradiction.

If $a_1a_2+a_2b_1+2a_1b_2=0$, then $b_2=-\frac{a_1a_2+a_2b_1}{2a_1}$. Substituting this into \eqref{di1} yields 
\begin{equation}
 8\frac{a_1}{a_2}(a_2-b_2)b_2^2+(3a_2-2b_2)b_2^2 X'=0,
\end{equation} 
which again implies $b_2=0$, a contradiction. The final case $2a_2+b_2=0$ leads to $b_1=3 a_1$, which in \eqref{di1} results 
$ 7 a_2 X'+24 a_1 =0$, a contradiction. 

If $b_2=a_1+b_1=0$, then $C_2=0$ becomes
\begin{equation*}
a_{2}^{3}(X')^{2}X+( \lambda a_{2}^{2}) (X')^{2} -4a_{2}^{2}b_{1} X'X
-4\lambda a_{2}b_{1} X'+( 4a_{2}b_{1}^{2}) X=0
\end{equation*}
Substituting $X$ into $C_1=0$ gives $2(a_{2}X'-2b_{1}) ^{2}=0$, a contradiction by the Claim L.

Thus, there exists some neighborhood $V_2 \subset V_1$ on which $(a_{2}Y+2\lambda )Y''+a_{1}a_{2}Y\neq0$.\\
 
\item[Case 2] From \eqref{y2}, suppose now that $%
3Y''+2a_{2}Y'-a_{1}=0$ in $V_{2}.$ Isolating $%
Y''$, we find 
\begin{equation}
Y''=-\frac{2a_{2}}{3}Y'+\frac{a_{1}}{3}
\label{y'' expression 2}.
\end{equation}%
From (\ref{xy}), we obtain the following expression for $Y'$
almost everywhere in $V_{2}:$ 
\begin{equation}\label{y5}
Y'=\frac{(a_{1}-3b_{1})Y}{(2a_{2}+3b_{2})Y+3\lambda }.
\end{equation}

As before, substituting into \eqref{eq3}, we find $\sum_{n=0}^3C_nY^n=0$, where the $C_n$ are all polynomials with respect to $X$ and 
$X'$. All of them must vanish identically. A computation of $C_0$ gives $C_0= -6a_{1}\lambda ^{2}X^{2}X'$. This implies $a_1=0$. In consequence, and from the expression \eqref{y5} of $Y'$, we deduce that $b_1\not=0$. The arguments follow similar steps as in the previous Case 1. From the condition $C_3=0$, we have the following equality:
\begin{equation}\label{y6}
 0=(2a_2+3b_2)(3a_2b_2XX'-6\lambda b_1-4a_2b_1X+2\lambda a_2X'+3\lambda b_2X'+2a_{2}^{2}XX').
\end{equation}
If $2a_2+3b_2=0$, then the conditions $C_2=0$ and $C_1=0$ lead to $b_2X-2\lambda=0$, which is a contradiction as $X$ would be constant.

If $2a_2+3b_2\neq0$, we can sove $X$ from \eqref{y6} and substitute it into $C_1=0$. After simplification, we obtain 
$$54b_1^2-3b_1(10a_2+21b_2)X'+(2a_2^2+15a_2b_2+18b_2^2)X'^2=0.$$
This implies $b_1=0$ which, together with the fact that $a_1=0$, implies $Y'=0$ in \eqref{y5}. This is a contradiction by Claim L. 

Thus, there exists some neighborhood $V_3 \subset V_2$ on which $3Y''+2a_{2}Y'-a_{1}\neq0$.

\item[Case 3] From \eqref{y3}, suppose $2\lambda Y''+\lambda a_{2}Y'+a_{1}a_{2}Y=0$
in $V_{3}$. This implies, 
\begin{equation*}
Y''=-\frac{a_{2}}{2}Y'-\frac{a_{1}a_{2}}{2\lambda }Y.
\end{equation*}
In particular, $a_2\not=0$. Using this in \eqref{xy}, we have almost everywhere in $V_{3}$:

\begin{equation*}
Y'=\frac{-a_{1}a_{2}Y^{2}-2\lambda b_{1}Y}{\lambda
((a_{2}+2b_{2})Y+2\lambda )},
\end{equation*}
and consequently,
\begin{equation}\label{y7}
Y''=-\frac{a_2Y((a_1-b_1)\lambda+a_1b_2Y)}{\lambda(2\lambda+(a_2+2b_2)Y)}.
\end{equation}
Substituting both equalities into (\ref{eq3}), we obtain a polynomial equation $\sum_{n=0}^3C_nY^n=0$, where the $C_n$ are polynomials in $X$ and $X'$. Thus, all $C_n$ must vanish. The coefficient $C_0$ is 
$$C_{0}=4 \lambda ^3 (a_1 X+\lambda X' ) (X' (a_2 X+\lambda )+b_1 X).$$
 \begin{enumerate}
 \item Subcase $a_1 X+\lambda X'=0$ identically. Then $X'=-\frac{a_1}{\lambda }X$. Substituting this into the other coefficients, we get 
\begin{equation*}
\begin{split}
C_1&=2a_2\lambda^2(a_1+b_1)^2X^2,\\
 C_2&=4a_1a_2(a_1+b_1)(a_2+b_2)\lambda X^2,\\
 C_3&=2a_1^2a_2(a_2+b_2)^2X^2.
 \end{split}
 \end{equation*}
From the expression of $Y''$ in \eqref{y7}, we have $a_2\not=0$, and $a_1$ and $b_1$ cannot both be zero. Thus we deduce $b_1=-a_1$ and $b_2=-a_2$. On the other hand, solving $a_1 X+\lambda X'=0$ yields $X(u)=m_1e^{-\frac{a_1}{\lambda}u}$. From the expression of $X''$ in \eqref{xy}, we get $a_2=-\frac{a_1}{\lambda}$. Using this information for the function $Y$, we obtain
 $Y'=\frac{a_1}{\lambda}Y$, which gives $Y(v)=m_2e^{\frac{a_1}{\lambda}v}$. This implies $YX'+XY'=0$, which is a contradiction.
 \item Subcase $a_1 X+\lambda X'\not=0$ and $X' (a_2 X+\lambda )+b_1 X=0$ identically. The arguments are similar as in the previous item and we omit the details. We have $X'=-\frac{b_1}{\lambda+a_2X}X$. In particular $b_1\not=0$. Substituting in $C_1$, $C_2$ and $C_3$, we obtain polynomial equations on $X$, whose coefficients must vanish. Hence, we find $b_1=a_1$, $b_2=0$. Finally, $Y''=0$ in \eqref{y7}, which is a contradiction by the Claim L.
\end{enumerate}

Therefore, there must exist some neighborhood $V_{4}\subset V_{3}$ on which $%
2\lambda Y''+\lambda a_{2}Y'+a_{1}a_{2}Y\neq 0$.

\item[Case 4] From \eqref{y4}, suppose $a_{2}Y'-2a_{1}=0$ in $V_{4}$. This implies that $Y$ is 
linear, a case already discarded. 
\end{description}

Consequently, there exists a neighborhood $V_{5}\subset V_{4}$ on which $a_{2}Y'-2a_{1}\neq 0$. In this neighborhood, the polynomial $p$ must be irreducible by construction. Thus, the polynomial
\begin{equation}\label{y8}
\begin{split}
q(X,X'):&=\frac{p(X,X')}{p_{6}} =(\frac{p_{3}}{p_{6}}X+1)\Big(X'^2+(%
\frac{p_{1}}{p_{6}}X^{2}+\frac{p_{4}}{p_{6}}X+\frac{p_{7}}{p_{6}})X'
\\
&+\frac{p_{2}}{p_{6}}X^{2}+\frac{p_{5}}{p_{6}}X\Big)
\end{split}
\end{equation}%
is also irreducible for all $v\in V_{5}$.

Taking two values $v_1,v_2\in V_5$, we obtain two irreducible polynomials $q_1,q_2$ of the same degree with the same root $(X,X')$ for all $u \in U_{0}$. By Lemma \ref{fulton}, either $(X,X')$ assumes finitely many values (which implies $X$ is constant, a contradiction) or there exists a function $\phi$ such that $q_1=\phi q_2$. Since both polynomials have the same constant term with respect to $X'$, it follows that $\phi = 1$ and thus $q_1=q_2$. By the arbitrariness of $v_1$ and $v_2$ we conclude that all coefficients of $q$ in \eqref{y8}are constant for all $v\in V_5$. That is, by this Lemma, either $X$ and $X'$ are constants in $U_{0}$, which
is a contradiction, or all of the coefficients $\frac{p_{j}}{p_{6}}$ are
constants in $V_{5}$. In particular, $\frac{p_{7}}{p_{6}}=2Y'$ must be constant in $V_{5}$, which implies $Y''=0$, a contradiction  by   Claim L. Therefore, we conclude that the
case $\lambda \neq 0$ is impossible.

This completes the proof of Theorem \ref{t1} for the case $H^2=1$.



\section*{Acknowledgements}
Rafael Belli would like to thank the Department of Geometry and Topology of the University of Granada for their hospitality, where this work was carried out. 
Rafael L\'opez has been partially supported by MINECO/MICINN/FEDER grant
no. PID2023-150727NB-I00, and by the ``Mar\'{\i}a de Maeztu'' Excellence
Unit IMAG, reference CEX2020-001105- M, funded by MCINN/AEI/10.13039/
501100011033/ CEX2020-001105-M.

\section*{Conflict of interest}

The authors have no relevant financial or non-financial interests to disclose.


\end{document}